\documentclass{conm-p-l}
\usepackage{amsmath}
\usepackage{amsthm}
\usepackage{amsfonts}
\usepackage{amssymb}

\newcommand{\1}{{1 \hspace{-0.35em} {\rm 1}}}
\newcommand{\I}{\mathcal{I}}

\newcommand{\OO}{\mathcal{O}}

\newcommand{\B}{\mathcal{B}}

\newcommand{\lan}{\langle}
\newcommand{\ra}{\rangle}
\newcommand{\End}{{\rm End}}
\newcommand{\Hom}{{\rm Hom}}
\newcommand{\ve}{\varepsilon}
\newcommand{\T}{\mathcal{T}}

\newcommand{\F}{\mathcal{F}}

\newcommand{\one}{{\rm Id}}

\newcommand{\Z}{\mathbb{Z}}
\newcommand{\g}{\mathfrak{g}}
\newcommand{\N}{\mathbb{N}}

\newcommand{\C}{\mathbb{C}}

\newcommand{\la}{{\lambda}}

\newcommand{\al}{\alpha}
\newtheorem{theorem}{Theorem}[section]
\newtheorem{lemma}[theorem]{Lemma}

\theoremstyle{definition}
\newtheorem{definition}[theorem]{Definition}
\newtheorem{example}[theorem]{Example}

\theoremstyle{remark}
\newtheorem{remark}[theorem]{Remark}

\numberwithin{equation}{section}

\begin{document}

\title{From Quantum Groups to Unitary Modular Tensor Categories}

\author{Eric C. Rowell}
\address{Department of Mathematics, Indiana University,
Bloomington, IN 47405}

\email{errowell@indiana.edu}

\date{January 31, 2005.}

\begin{abstract} Modular tensor categories are generalizations of the representation categories of
quantum groups at roots of unity axiomatizing the properties
necessary to produce 3-dimensional TQFTs.  Although other
constructions have since been found, quantum groups remain the most
prolific source.  Recently proposed applications to quantum
computing have provided an impetus to understand and describe these
examples as explicitly as possible, especially those that are
``physically feasible."  We survey the current status of the problem
of producing unitary modular tensor categories from quantum groups,
emphasizing explicit computations.
\end{abstract}

\maketitle

\section{Introduction} We outline the development of the theory of
modular tensor categories from quantum groups with an eye towards
new applications to quantum computing that motivate our point of
view.
  In this article, we take
\emph{quantum group} to mean the ``classical" $q$-deformation of the
universal enveloping algebra of a simple complex finite dimensional Lie algebra as in the book
by Lusztig \cite{Lusztigbook}, rather than the broader class of Hopf
algebras the term sometimes describes.

\subsection{Background}
The representation theory of quantum groups has proven to be a
useful tool and a fruitful source of examples in many areas of
mathematics.  The general definition of a quantum group (as a Hopf
algebra) was given around 1985 by Drinfeld \cite{D} and
independently Jimbo \cite{Jimbo} as a method for finding solutions
to the quantum Yang-Baxter equation.  These solutions led to new
representations of Artin's braid group $\B_n$ and connections to
link invariants.  In fact, specializations of the famous polynomial
invariants of Jones \cite{Jones85}, the six-authored paper
\cite{HOMFLY} and Kauffman \cite{Kauf} have been obtained in this
way.  Reshetikhin and Turaev \cite{ReshTur} used this connection to
derive invariants of 3-manifolds from \emph{modular Hopf algebras},
examples of which can be found among quantum groups at roots of
unity (see \cite{ReshTur} and \cite{TW93}, much simplified by
constructions in \cite{andersen}).  When Witten \cite{Witten88}
introduced the notion of a topological quantum field theory (TQFT)
relating ideas from quantum field theory to manifold invariants,
non-trivial examples were immediately available from the
constructions in \cite{ReshTur} (after reconciling  notation).
Modular Hopf algebras were replaced by the more general framework of
modular tensor categories (MTCs) by Turaev \cite{Tur92} (building on
definitions in \cite{MacLane} and \cite{JS}), axiomatizing the
conditions necessary (and sufficient, see \cite{TuraevBook},
Introduction) to construct 3-dimensional TQFTs.

Aside from the quantum group approach to MTCs, there are several
other general constructions.  Representation categories of Hopf
algebra doubles of finite group algebras are examples of MTCs that
are often included with quantum groups in the more general
discussion of Hopf algebra constructions. A geometric construction
using link invariants and tangle categories was introduced in
\cite{TuraevBook},
 advanced by Turaev and Wenzl in \cite{TW97} and somewhat simplified by
Blanchet and Beliakova in \cite{BB}. However, all examples that
have been carried out lead to MTCs also obtainable from quantum
groups.
 Yet another construction of MTCs from representations of
 vertex operator algebras has
recently appeared in a paper by Huang \cite{Huang}. See Subsection \ref{other} for
further discussion of these approaches.

 Although it is expected that there are non-trivial examples of MTCs
that do not arise from Hopf algebras (\emph{e. ~g.} quantum groups
and finite group algebras), none have been rigorously produced. This is
probably due to the highly advanced state of the theory of
representations of quantum groups at roots of unity provided by the
pioneering work of many including Lusztig (\cite{Lusztigbook}) and
Andersen and his co-authors (\cite{andersen}, \cite{AndPar}
\cite{APW}).  The description of the MTCs derived from quantum
groups can be understood with little more than a firm grasp on the
theory of representations of simple finite-dimensional Lie algebras
found in Humphrey's book \cite{Hum} or any other introductory text.

The purpose of this paper is two-fold: to survey what is known about
the modularity and unitarity of categories arising from quantum
groups at roots of unity, and to give useful combinatorial tools for
explicit computations in these categories.  For more in-depth
developments the reader is directed to two references: 1) Bakalov
and Kirillov's text (\cite{BK}, Sections 1.3 and 3.3)
 contains concise constructions and examples of quantum group MTCs, and 2) Sawin's paper
 (\cite{S2}) gives a thorough treatment of the representation theoretic details
 necessary to construct MTCs from quantum groups.  The
 modularity results described below partially overlap with Section 6 of
 \cite{S2}.
\subsection{Motivation}
There are two fairly well-known motivations for studying MTCs.  They
are applications to low-dimensional topology (see
\cite{TuraevBook}), and conformal field theory (see \cite{Huang} and
references therein).

Recently, an application of unitary MTCs to quantum computing has
been proposed by Freedman and Kitaev and advanced in the series of
papers (\cite{FKW}, \cite{FKLW}, \cite{FLW} and \cite{FNSWW}). Their
topological model for quantum computing has a major advantage over
the ``classical" qubit model in that errors are corrected on the
physical level and so has a higher error threshold.  For a very
readable introduction to topological quantum computing see
\cite{FKW}.   In this model unitary MTCs play the role of the
software, while the hardware is implemented via a quantum physical
system and the interface between them is achieved by a 3-dimensional
TQFT.  The MTCs encode the symmetries of the corresponding physical
systems (called topological states, see \cite{FNSWW}), and must be
unitary by physical considerations.

Aside from the problem of constructing unitary MTCs, there are
several open problems currently being studied related to the quantum
computing applications.  One question is whether the images of the
irreducible unitary braid representations (see Remark
\ref{braidrem}) afforded by a unitary MTC are dense in the unitary
group.  This is related to a \emph{sine qua non} of quantum
computation known as \emph{universality}.  Progress towards
answering this question has been made in \cite{FLW} and was extended
by Larsen, Wang, and the author in \cite{LRW}.  Another problem is
to prove the conjecture of Z. Wang: \emph{There are finitely many
MTCs of a fixed rank} (see Subsection \ref{notrems}). This has been
verified for ranks 1,2,3 and 4: see \cite{Ostrik1} and
\cite{Ostrik2} for ranks 2 and 3 respectively, and \cite{BRSW} for
both ranks 3 and 4. It is with this conjecture in mind that we
provide generating functions for ranks of categories in Subsection
\ref{gen}.

\subsection*{Acknowledgements} The author wishes to thank the referees
for especially careful readings of previous versions of this article
and for comments leading to a much-improved
exposition.  Special thanks
also to Z. Wang for many useful discussions on
topological quantum computation.

\section{General Definitions}
We give the basic categorical definitions for modular tensor
categories, remark on some consequences and describe the crucial
condition of unitarity.
\subsection{Axioms}\label{cat}
In this subsection we outline the axioms for the categories we are
interested in. We follow the paper \cite{Tur92}, and refer to that
paper or the books
by Turaev \cite{TuraevBook} or Kassel \cite{K} for a complete treatment.\\

Let $\OO$ be a category defined over a subfield $k\subset\C$. A
modular tensor category is a semisimple ribbon \emph{Ab}-category
$\OO$ with finitely many isomorphism classes of simple objects
satisfying a non-degeneracy condition. We unravel these adjectives
with the following definitions.
\begin{enumerate}
\item A \textbf{monoidal category} is a category with
 a tensor product $\otimes$
and an identity object $\1$ satisfying axioms that guarantee that
the tensor product is associative (at least up to isomorphism) and
that
\begin{equation*}\label{idob}
\1\otimes X\cong X\otimes\1\cong X
\end{equation*}
 for any object $X$.  See \cite{MacLane} for details.
\item A monoidal category has \textbf{duality} if there is a dual
object $X^*$ for each object $X$ and morphisms
$$b_X: \1 \rightarrow X\otimes X^*, d_X: X^*\otimes X \rightarrow \1$$
satisfying
\begin{eqnarray*}
(\one_X\otimes d_X)(b_X\otimes \one_X)&=&\one_X, \\
(d_X\otimes \one_{X^*})(\one_{X*}\otimes b_X)&=&\one_{X^*}.
\end{eqnarray*}  The duality
allows us to define duals of morphisms too: for any\\ $f\in
\Hom(X,Y)$ we define $f^*\in \Hom(Y^*,X^*)$ by:
$$f^*=
(d_Y\otimes \one_{X^*})(\one_{Y^*}\otimes f\otimes
\one_{X^*})(\one_{Y^*}\otimes b_X).$$
\item A
\textbf{braiding} in a monoidal category is a natural family of isomorphisms
$$c_{X,Y}: X\otimes Y \rightarrow Y\otimes X$$
satisfying
\begin{eqnarray*}
c_{X,Y\otimes Z}&=&(\one_Y\otimes c_{X,Z})(c_{X,Y}\otimes \one_Z), \\
c_{X\otimes Y,Z}&=&(c_{X,Z}\otimes \one_Y)(\one_X\otimes c_{Y,Z}).
\end{eqnarray*}
\item A \textbf{twist} in a braided monoidal category is a natural family of isomorphisms
$$\theta_X: X \rightarrow X$$
satisfying:
\begin{eqnarray*}
\theta_{X\otimes Y}&=&c_{Y,X}c_{X,Y}(\theta_X\otimes \theta_Y).
\end{eqnarray*}
\item In the presence of a braiding, a twist and duality these structures
are compatible if
\begin{eqnarray*}
\theta_{X^*}=(\theta_X)^*.
\end{eqnarray*}  A braided monoidal category with a twist and a compatible duality is
a \textbf{ribbon} category.
\item An \textbf{$Ab$-category} is a monoidal category in which all morphism spaces
are $k$-vector spaces and the composition and tensor product of
morphisms are bilinear.
\item An $Ab$-category is \textbf{semisimple} if it has the property that every object
$X$ is isomorphic to a finite direct sum of \emph{simple}
objects--that is, objects $X_i$ with $\End(X_i)\cong k$ satisfying the conclusion of Schur's Lemma:
$$\Hom(X_i,X_j)=0\quad \text{for $i\neq j$}.$$  Turaev \cite{TuraevBook}
gives a weaker condition for semisimplicity avoiding
direct sums, but we omit it for brevity's sake.
\item In a ribbon \emph{Ab}-category one may define a $k$-linear \textbf{trace} of
endomorphisms.  Let $f\in\End(X)$ for some object $X$.  Set:
\begin{eqnarray*}
tr(f)=d_Xc_{X,X^*}(\theta_Xf\otimes\one_{X^*})b_X
\end{eqnarray*}
where the right hand side is an element of $\End(\1)\cong k$. The
value of $tr(\one_X)$ is called the categorical dimension of $X$ and
denoted $\dim(X)$.
\item A semisimple ribbon \emph{Ab}-category is called a \textbf{modular tensor category}
if it has finitely many isomorphism classes of simple objects
enumerated as $\{X_0=\1,X_1,\ldots,X_{n-1}\}$ and the so called
\emph{$S$-matrix} with entries $$S_{i,j}:=tr(c_{X_j,X_i}\circ
c_{X_i,X_j})$$ is invertible. Observe that $S$ is a symmetric
matrix.
\end{enumerate}

\subsection{Notation and Remarks}\label{notrems}
In a semisimple ribbon \emph{Ab}-category $\OO$ with finitely many
simple classes the set of simple classes generates a semiring over
$k$ under $\otimes$ and $\oplus$.  This ring is called the
\emph{Grothendieck semiring} and denoted $Gr(\OO)$.  If
$\{X_0=\1,X_1,\ldots,X_{n-1}\}$ is a set of representatives of these
isomorphism classes, the \textbf{rank} of $\OO$ is $n$.  The axioms
guarantee that we have (using Kirillov's notation \cite{kirillov}):
\begin{eqnarray}\label{fusion}
X_i\otimes X_j\cong\sum_{k}N_{i,j}^kX_k
\end{eqnarray}
for some $N_{i,j}^k\in\N$.  These structure coefficients of
$Gr(\OO)$ are called the \emph{fusion coefficients} of $\OO$ and (\ref{fusion}) is sometimes called
a fusion rule.
  Having fixed
an ordering of simple objects as above, the fusion coefficients
give us a representation of $Gr(\OO)$ via $X_i\rightarrow N_i$ where
$N_i$, $(N_i)_{k,j}=(N_{i,j}^k)$ is called the \emph{fusion matrix} associated to
$X_i$.  If we denote by $i^*$ the index of the simple object
$X_i^*$, the braiding and associativity constraints give us:
\begin{eqnarray*}
N_{i,j}^k=N_{j,i}^k&=&N_{i,k^*}^{j^*}=N_{i^*,j^*}^{k^*},\\
N_{i,j}^0&=&\delta_{i,j^*}.
\end{eqnarray*}
It also follows from associativity that the fusion matrices pairwise
commute, so that full fusion rules may sometimes be computed just
from a single fusion matrix (\emph{i.e.} using a Gr\"obner basis
algorithm).

The first column (and row) of the $S$-matrix consists of the
categorical dimensions of the simple objects, \emph{i.e.}
$S_{i,0}=\dim(X_i)$.  We denote these dimensions by $d_i$.  We also
have that $S_{i,j}=S_{j,i}=S_{i^*,j^*}$. Since the twists
$\theta_X\in\End(X)$ for any object $X$, $\theta_{X_i}$ is a scalar
map (as $X_i$ is simple).  We denote this scalar by $\theta_i$.

Standard arguments show that the entries of the $S$-matrix are
determined by the categorical dimensions, the fusion rules and the
twists on these simple classes, giving the following
extremely useful formula (see \cite{BK}):
\begin{eqnarray}\label{sformula}
S_{i,j}=\frac{1}{\theta_i\theta_j}\sum_kN_{i^*,j}^kd_k\theta_k.
\end{eqnarray}
Provided $\OO$ is modular the $S$-matrix determines the fusion rules
via the \emph{Verlinde formula} (see \cite{BK}, and \cite{Huang}).
To express the formula we must introduce the quantity $D^2=\sum_i
d_i^2$. Then:
\begin{eqnarray}
N_{i,j}^k=\sum_t\frac{S_{i,t}S_{j,t}S_{k^*,t}}{D^2S_{0,t}}.
\end{eqnarray}
This formula corresponds to the following fact: the columns of the
$S$-matrix are simultaneous eigenvectors for the fusion matrices
$N_i$, and the categorical dimensions are eigenvalues.

\begin{remark}\label{braidrem}
The braiding morphisms $c_{X,X}$ induce a representation of the braid group $\B_n$
on $\End(X^{\otimes n})$ for any object $X$ via the operators
\begin{eqnarray*}R_i=
\one_X^{\otimes i-1}\otimes c_{X,X}\otimes\one_X^{\otimes n-i-1}\in\End(X^{\otimes n})\end{eqnarray*}
 and the generators $\sigma_i$ of $\B_n$ act by left composition by $R_i$.
\end{remark}
\begin{remark}
The term ``modular" comes from the following fact: if we set
$T=(\delta_{i,j}\theta_i)_{ij}$ then the map:
$$\begin{pmatrix}
0 &  -1\\
 1 &  0
\end{pmatrix}\rightarrow S, \begin{pmatrix}
1 &  1\\
 0 &  1
\end{pmatrix}\rightarrow T$$
defines a projective representation of the \emph{modular}
\emph{group} $SL(2,\Z)$.  In fact, by re-normalizing $S$ and $T$ one
gets an honest representation of $SL(2,\Z)$.
\end{remark}
\subsection{Unitarity}
A \textbf{Hermitian} ribbon \emph{Ab}-category has a conjugation:
$$\dag: \Hom(X,Y)\rightarrow \Hom(Y,X)$$
such that $(f^\dag)^\dag=f$, $(f\otimes g)^\dag=f^\dag\otimes
g^\dag$ and $(f\circ g)^\dag=g^\dag\circ f^\dag$.  On $k\subset\C$,
$\dag$ must also act as the usual conjugation.  Furthermore, $\dag$
must be compatible with the other structures present
\emph{i.e.}
\begin{eqnarray*}
(c_{X,Y})^\dag&=&(c_{X,Y})^{-1}, \\
(\theta_X)^\dag&=&(\theta_X)^{-1},\\
(b_X)^\dag&=&d_Xc_{X,X^*}(\theta_X\otimes \one_{X^*}),\\
(d_X)^\dag&=&(\one_{X^*}\otimes\theta_X^{-1})(c_{X^*,X})^{-1}b_X.
\end{eqnarray*}
For Hermitian ribbon \emph{Ab}-categories the categorical dimensions
$d_i$ are always real numbers. If in addition the Hermitian form
$(f,g)=tr(fg^\dag)$ is positive definite on $\Hom(X,Y)$ for any two
objects $X,Y\in\OO$, the category is called \textbf{unitary}, and
the categorical dimensions are positive real numbers.  If $\OO$ is
unitary, then the morphism spaces $\End(X)$ are Hilbert spaces with
the above form, and the representations
$$
\B_n\rightarrow\End(X^{\otimes n})$$ described in Remark \ref{braidrem} are unitary.

\section{Constructions}
MTCs have been derived in varying degrees of detail from several
sources.  A very general approach is through representations of
quantum groups at roots of unity.  We give a very broad outline of
how these are obtained and mention a few other sources and
constructions.
\subsection{MTCs from Quantum Groups}  The following construction is now standard,
and can be found in more detail in the books by Turaev
\cite{TuraevBook} or Bakalov and Kirillov Jr. ~\cite{BK} (both of which include examples).  The
 procedure is a culmination of the work of many, but the
major contributions following those of Drinfeld and Jimbo were from
Lusztig (see \cite{Lusztigbook}), Andersen and his collaborators
(\cite{APW},\cite{andersen} and \cite{AndPar}) and Turaev with
Reshetikhin (\cite{ReshTur}) and Wenzl (\cite{TW93}).
 Let $\g$ be a Lie algebra from one of the infinite families $ABCD$
or an exceptional Lie algebra of type $E,F$ or $G$ and $q$ a complex
number such that $q^2$ is a primitive $\ell$th root of unity, where
$\ell$ is greater than or equal to the dual Coxeter number of $\g$.
Let $U=U_q(\g)$ be Lusztig's \cite{Lusztigbook} ``modified form" of
the Drinfeld-Jimbo quantum group specialized at $q$ and denote by
$\T$ Andersen's \cite{andersen} subcategory of \emph{tilting
modules} over $U$.  A module $V$ is called \emph{tilting} if both
$V$ and its dual, $V^*$, admit Weyl filtrations: \emph{i.e.}
sequences $$ \{0\}=V_0\subset V_1\subset\cdots\subset V_n=V$$ with
each $V_i/V_{i-1}$ a Weyl module. The ratio of the square lengths of
a long root to a short root will play an important role in the
sequel, so we denote it by the letter $m$. Observe that $m=1$ for
Lie types $ADE$, $m=2$ for Lie types $BC$ and $F$, and $m=3$ for Lie
type $G$.  It can be shown that $\T$ is a (non-semisimple) ribbon
\emph{Ab}-category (see \cite{andersen} and \cite{TW93}). The ribbon
structure on $\T$ comes from the (ribbon) Hopf algebra structure on
$U$ (see \cite{ChP}), \emph{i.e.} the antipode, comultiplication,
$R$-matrix, quantum Casimir \emph{etc}. The set of indecomposable
tilting modules with $\dim(X)=0$ (categorical dimension) generates a
tensor ideal $\I\subset\T$, and semisimplicity is recovered by
taking the quotient category $\F=\T/\I$.  Moreover, the category
$\F$ has only finitely many isomorphism classes of simple objects,
labelled by the subset of dominant weights (denoted $P_+$) in the
\emph{fundamental alcove}:
$$C_\ell(\g):=\begin{cases}\{\la\in P_+: \lan\la+\rho,\vartheta_0\ra<\ell\} &
\text{if $m\mid \ell$}\\
\{\la\in P_+: \lan\la+\rho,\vartheta_1\ra<\ell\} & \text{if $m\nmid
\ell$}\end{cases}$$ where $\vartheta_0$ is the highest root and
$\vartheta_1$ is the highest short root. Here the form $\lan \: ,\,
\ra$ is normalized so that $\lan \alpha, \alpha \ra=2$ for
\emph{short} roots.  While $\F$ is always a semisimple Hermitian
ribbon \emph{Ab}-category with finitely many isomorphism classes of
simple objects, the further properties (modularity and unitarity) of
$\F$ depend on $\g$, the divisibility of $\ell$ by $m$, and the
specific choice of $q$.  We denote the category $\F$ by
$\mathcal{C}(\g,\ell,q)$ to emphasize this dependence.
 The $S$-matrices for these categories are well-known.  For $\la,\mu \in C_\ell(\g)$ we have:
\begin{eqnarray}\label{smat}
S_{\la,\mu}=\frac{\sum_{w\in W}\ve(w)q^{2\lan
\la+\rho,w(\mu+\rho)\ra}}{\sum_{w\in
W}\ve(w)q^{2\lan\rho,w(\rho)\ra}}
\end{eqnarray} where $\rho$ is the half sum of the positive roots and $\ve(w)$ denotes
the sign of the Weyl group element $w$.
\begin{remark}
In practice, Formula (\ref{sformula}) is often more
useful than Formula (\ref{smat}) for computing the entries of the $S$-matrix,
as computing the twists $\theta_\la$,
$q$-dimensions $d_\la$ (see below) and fusion coefficients $N_{\la,\mu}^\nu$
(via the quantum Racah formula, see \cite{AndPar} and \cite{S2})
is more straightforward than summing
over the Weyl group.
\end{remark}
The twist coefficients for simple objects are also well known:
$\theta_\la=q^{\lan \la,\la+2\rho\ra}$, as are the categorical
$q$-dimensions:

$$d_\la=\prod_{\al\in\Phi_+}\frac{[\lan\la+\rho,\al\ra]}{[\lan\rho,\al\ra]}$$ where
$[n]=\frac{q^n-q^{-n}}{q-q^{-1}}$ and $\Phi_+$ is the set of
positive roots.

We note that the fusion coefficients of $\mathcal{C}(\g,\ell,q)$ only depend on $\g$ and
$\ell$. A
complete description of the braiding and associativity maps is quite
difficult in general; fortunately one is usually content to know
they exist, relying on the $S$-matrix, fusion matrices and twists
for most calculations.
\begin{remark} An issue has recently come to light regarding the
explicit fusion rules for these categories.  While
Andersen-Paradowski \cite{AndPar} showed that for many cases the
fusion rules for the truncated tensor product in the category $\F$
are determined from the classical multiplicities by an
anti-symmetrization over the affine Weyl group, their proof appeared
in a paper that restricted attention to the root lattice.  Evidently
the first
 general proof of this ``quantum Racah" formula is in the preprint \cite{S2}.
\end{remark}

\subsection{Other Constructions}\label{other}
The most direct construction of MTCs comes from the representation
category of the semidirect product $D(G):=k[G]\ltimes \F(G)$ of the
group algebra of a finite group with its (Hopf algebra) dual and can
be found in the book \cite{BK}.  For example, the representation
category of the Hopf algebra $D(S_3)$ is a rank 8 MTC that does not
arise from a quantum group construction as outlined above.  These
MTCs always have integer $q$-dimensions.

  The geometric construction of MTCs alluded to in the introduction
is summarized as follows. One starts with a link invariant
satisfying a number of mild (but technical) conditions and produces
a new category from the category of tangles via an \emph{idempotent
completion} of quotients of endomorphism spaces.  This produces a
semisimple braided category, and if there is explicit information
available for the link invariant one can sometimes verify the
remaining axioms. This has been carried out for the Jones polynomial
(Chapter XII of \cite{TuraevBook}) and the Kauffman polynomial
\cite{TW97}. Blanchet and Beliakova \cite{BB} gave a complete
analysis of the modularity and modularizability of these categories
corresponding to $BMW$ algebras--the algebras supporting the
Kauffman polynomial. Although the work in \cite{BB} eliminated the
need to appeal to quantum group characters as in \cite{TW97}, these
constructions give rise to essentially the same MTCs as those
obtained from quantum groups of types $B$, $C$ and $D$ at roots of
unity.  An advantage of this geometric approach is that the braid
representations are more transparent than in the quantum group
construction, although one pays for this convenience by having a
less natural description of objects.

 As we noted in the introduction, MTCs
have also been constructed from representation categories of
certain vertex operator algebras (VOAs) by
Huang \cite{Huang}. Rigidity and modularity are the most difficult
to verify, while the monoidal structure was previously obtained.  The allure of this
approach is that it
includes a proof of a very general form of the Verlinde conjecture from conformal
field theory.  Although this VOA construction of MTCs is more difficult than other
approaches, it gives credence to the thesis that MTCs describe symmetries in quantum physical
systems.

There are two indirect constructions that should be mentioned.  One
is the quantum double technique of M\"uger \cite{Muger} (inspired by
the double of a Hopf algebra) by which an MTC is
constructed by ``doubling" a monoidal category with some further
technical properties.  An example of this approach is the finite
group algebra construction mentioned above.  Brugui\`eres \cite{Br}
describes conditions under which one may \emph{modularize} a
category that satisfies all of the axioms of an MTC except the
invertibility of the $S$-matrix (called a \textbf{pre-modular category}).  This
corresponds essentially to taking a quotient or sub-category that
does satisfy the modularity axiom.

\section{Modularity, Unitarity and Ranks for Quantum Groups}
There remains a fair amount of work to be done to have a complete
theory of abstract unitary modular tensor categories; however, for
quantum groups much is known.  The condition of modularity has been
settled for nearly all of the categories $\mathcal{C}(\g,\ell,q)$,
as well as the question of unitarity.

The modularity condition is often difficult to verify.  Recently a
\textbf{modularity criterion} was proved that sometimes simplifies
the work (see \cite{Br}):
\begin{theorem}[Brugu\`ieres]\label{brug}
Suppose $\OO$ is a pre-modular category, and
let $\{X_0=\1,X_1,\ldots,X_{n-1}\}$ be a set of representatives of the
simple isomorphism classes.  Then $\OO$ is modular if and only if
 $$\mathcal{N}:=\{X_i: S_{i,j}=d_{i}d_{j} \text{ for all $X_j$}\}=\{\1\}.$$
\end{theorem}
Observe that one has $S_{0,j}=d_0d_j=d_j$.  The non-trivial elements of
$\mathcal{N}$ are the \textbf{obstructions to modularity}, \emph{i.e.} the objects
for which the corresponding columns in the $S$-matrix are scalar
multiples of the first column.

In the following subsections we describe the modularity and
unitarity of the categories $\mathcal{C}(\g,\ell,q)$, first for the
cases that can be handled uniformly, and then for those that must be
considered individually as well as a few subcategories of interest.
Subsection \ref{gen} concerns the ranks of the categories
$\mathcal{C}(\g,\ell,q)$ and can be safely skipped by those readers
not interested in this issue.

\subsection{Uniform Cases $m \mid \ell$} For the categories $\mathcal{C}(\g,\ell,q)$
the cases where $\ell$ is divisible by $m$ have been mainly studied in the literature.
The invertibility of $S$ for Lie types $A$ and $C$ with $q=e^{\pi i/ \ell}$ was shown
 in \cite{TW93} using the work of Kac and Peterson \cite{KP}, and a complete treatment
 (for all Lie types with $q=e^{\pi i/ \ell}$) is found in \cite{kirillov}.
The invertibility can be extended to other values of $q$ by the
following Galois argument, which is found in \cite{TW97} in a
different form. By Formula (\ref{smat}) we see that the entries of
the $S$-matrix:
$$S_{\la\mu}=(const.)\sum_{w\in W}\ve(w)q^{2\lan w(\la+\rho),\mu+\rho\ra}$$
are polynomials in $q^{1/d}$ where $d\in\N$ is minimal so that
$d\lan\la,\mu\ra\in\Z$ for all weights $\la,\mu$.  Thus $\det(S)$ is
non-zero for any Galois conjugate of $e^{\pi i/d\ell}$, \emph{i.e.}
for any $q=e^{z\pi i/\ell}$ with $\gcd(z,d\ell)=1$. Table
\ref{dvalues} lists the values of $d$ for all Lie types for which
$d\not=1$. Notice that there are sub-cases for types $B$ and $D$.
\begin{table}
\begin{tabular}{{|c||}*{6}{c|}}
\hline
$X_r$ & $A_r$ & $B_r$, $r$ odd &  $D_r$, $r$ even& $D_r$, $r$ odd &$E_6$& $E_7$\\
\hline $d$& $r+1$&  $2$&  $2$&  $4$& $3$& $2$\\
 \hline
\end{tabular} \caption{Values of $d$ for Lie algebra types with $d\not=1$}\label{dvalues}
\end{table}
When $d=1$ and $m|\ell$ the uniform case covers all possibilities,
since then the condition $\gcd(z,d\ell)=1$ is equivalent to the
original assumption that $q^2$ is a primitive $\ell$th root of
unity. So when $m|\ell$, the cases $B_r$ with $r$ even, $C_r$,
$E_8$, $F_4$, and $G_2$ do not require further attention.  If $m=1$
and $\gcd(\ell,d)\not=1$ the condition $\gcd(z,d\ell)=1$ also
degenerates to the original assumption that $q^2$ is a primitive
$\ell$th root of unity so we need not consider $D_r$ with $\ell$
even, $E_6$ with $3|\ell$ or $E_7$ with $\ell$ even.

Following a conjecture of Kirillov Jr. ~\cite{kirillov}, Wenzl \cite{WenzlCstar}
showed that the Hermitian form on
$\mathcal{C}(\g,\ell,q)$ is positive definite for the uniform cases
for certain values of $q$, and Xu \cite{Xu} independently showed
some of the cases covered by Wenzl.  Their results are summarized
in:
\begin{theorem}[Wenzl/Xu]\label{cstar}
The categories $\mathcal{C}(\g,\ell,q)$ are unitary when $m|\ell$
and $q=e^{\pi i/\ell}$.
\end{theorem}
\subsection{Type $A$}
For Lie type $A_r$ corresponding to $\g=\mathfrak{sl}_{r+1}$ we have
$m=1$ and $d=r+1$.  Brugui\`eres \cite{Br} shows that one has
modularity for $q=e^{z\pi i/\ell}$ \emph{if and only if}
$\gcd(z,(r+1)\ell)=1$. Moreover, Masbaum and Wenzl \cite{MaW} show
that when $\gcd(\ell,r+1)=1$ the subcategory of
$\mathcal{C}(\mathfrak{sl}_{r+1},\ell,q)$ generated by the simple
objects labelled by integer weights is a modular subcategory whose
rank is $1/(r+1)$ times the rank of the full category.  There are a
number of other proofs of this fact, see e.g. \cite{Br} Section 5.
Denote this subcategory by $\Z(A_r)$, and see Subsection \ref{ex1}.

\subsection{Type $B$, $\ell$ odd}
The category $\mathcal{C}(\mathfrak{so}_{2r+1},\ell,q)$ with $\ell$
odd has been considered to some extent by several authors including
Sawin \cite{Sawin}, \cite{S2} and Le-Turaev \cite{LeTuraev}.  It is
shown in (\cite{TW97}, Theorem 9.9) that if $\ell$ is odd, the
subcategory of $\mathcal{C}(\mathfrak{so}_{2r+1},\ell,q)$ generated
by simple objects labelled by integer weights is modular and has
rank exactly half of that of
$\mathcal{C}(\mathfrak{so}_{2r+1},\ell,q)$.  Combining the
computations in \cite{Ribbon} and the modularity criterion of
\cite{Br} one has:
\begin{theorem}\label{bmod}
The category $\mathcal{C}(\mathfrak{so}_{2r+1},\ell,q)$ with $\ell$
odd is modular if and only if $q^\ell=-1$ and $r$ is odd.
\end{theorem}
\begin{proof}
By the modularity criterion we wish to show that there are
obstructions to modularity (\emph{i.e.} non-trivial objects in the
set $\mathcal{N}$, see Theorem \ref{brug}) if and only if the
conditions of the theorem are not satisfied. By the modularity of
the subcategory generated by simple objects labelled by integer
weights, any obstructing object must be labelled by a half-integer
weight. In \cite{Ribbon} the object $X_\gamma$ labelled by the
(half-integer) weight that is furthest from the $0$ weight in the
fundamental alcove is shown to induce an involution of the
fundamental alcove (by tensoring with $X_\gamma$) that preserves
$q$-dimension up to a sign.  This implies that $X_\gamma$ is the
only potential obstruction to modularity.  In \cite{Ribbon}
(Scholium 4.11) the signs of the $q$-dimensions are analyzed, and
 the theorem then follows from the explicit
computations of $N_{\gamma,\la}^\nu$, $d_\la$ and $\theta_\la$ (also found in \cite{Ribbon})
together with Formula (\ref{sformula}) and the obstruction
equation $S_{\gamma,\la}=d_\gamma d_\la$.
\end{proof}

The subject of the author's thesis \cite{thesis} (the results of
which can be found in \cite{Ribbon}) is the question of
unitarizability of the family of categories
$\mathcal{C}(\mathfrak{so}_{2r+1},\ell,q)$ with $\ell$ odd.  Using
an analysis of the characters of the Grothendieck semirings it is
shown that no member of this family of categories is unitary.  In
fact, there is a much stronger statement, for which we need the
following definition:
\begin{definition}
A pre-modular category $\OO$ is called \textbf{unitarizable} if
$\OO$ is tensor equivalent to a unitary pre-modular category
$\OO^\prime$.  By tensor equivalent we mean there exists a functor preserving
the monoidal structure that is
bijective on morphisms and such that every object in the target category is
isomorphic to an object in the image of the functor.
\end{definition}

Using a structure theorem of Tuba and Wenzl \cite{TubaWenzl} it is
shown in \cite{thesis} that:
\begin{theorem}\label{bunit}
Fix $q$ with $q^2$ a primitive $\ell$th root of unity, $\ell$ odd,
and $r$ satisfying $2(2r+1)<\ell$. Then \emph{no} braided tensor
category having the same Grothendieck semiring as
$\mathcal{C}(\mathfrak{so}_{2r+1},\ell,q)$ is unitarizable.
\end{theorem}
\begin{remark} When $\ell<2(2r+1)$ the rank of $\mathcal{C}(\mathfrak{so}_{2r+1},\ell,q)$ is
relatively small and the fusion rules of the category may coincide
with those of another category that is known to be unitarizable. For
example $\mathcal{C}(\mathfrak{so}_{5},7,q)$ has the same
Grothendieck semiring as $\mathcal{C}(\mathfrak{sl}_{2},7,q)$ which
is unitary for $q=e^{\pi i/7}$.
\end{remark}

\subsection{Type $C$, $\ell$ odd}
For type $C$ one has $m=2$, so it remains to analyze the cases with
$\ell$ odd. For this, we resort to the ``rank-level duality" result
of \cite{Ribbon} (Corollary 6.6) showing that the categories
$\mathcal{C}(\mathfrak{so}_{2r+1},\ell,q)$ and
$\mathcal{C}(\mathfrak{sp}_{\ell-2r-1},\ell,q)$ are tensor
equivalent. Theorem \ref{bunit} immediately implies these categories
are not unitarizable for $\ell$ odd if $2(2r+1)<\ell$. Moreover, the
technique in the proof of Theorem \ref{bmod} can be applied to this
case using the explicit values of $d_\la$ and $\theta_\la$ and the
image of the
 object $X_\gamma$ under the tensor equivalence afforded
by this rank-level duality.  We
then have:
\begin{theorem}
If $\ell$ is odd, the categories
$\mathcal{C}(\mathfrak{sp}_{2r},\ell,q)$ are not modular and if in
addition $2(2r+1)<\ell$ they are not unitarizable.
\end{theorem}

\subsection{Remaining Types $D$, $E_6$ and $E_7$ Cases}
The only remaining question for the sub-cases not covered
by the uniform case is whether the condition $\gcd(z,d\ell)=1$ is
necessary for modularity.  For Lie types $D$ and $E_7$ the sub-cases
correspond to $\ell$ odd, and for Lie type $E_6$ the sub-cases correspond to $3\nmid \ell$.
In our opinion this question is still open, of limited
interest and one probably does not get modularity.
\subsection{Types $F_4$ with $\ell$ odd, and $G_2$ with $3\nmid\ell$}
To our knowledge both the question of modularity and
unitarizability are still open for $F_4$ with $\ell$ odd and $G_2$
with $3\nmid\ell$.  In light of the results in the Lie types $B$ for $\ell$ odd
(see Theorems \ref{bmod} and \ref{bunit}),
 one might expect to find that these categories
are not unitarizable (except possibly for small $\ell$), but
sometimes modular.

\subsection{Generating Functions for $|C_\ell(\g)|$}\label{gen}
For
applications it is useful to know the ranks of the categories
$\mathcal{C}(\g,\ell,q)$.

We define an
auxiliary label $\ell_m=0$ if $m\mid\ell$ and $\ell_m=1$ if
$m\nmid\ell$ for notational convenience.
We reduce the problem of determining the cardinalities of the
labeling sets $C_\ell(\g)$ of simple objects to counting partitions of $n$ with parts
in a fixed (finite) multiset $\mathcal{S}(\g,\ell_m)$ that depends
only on the rank and Lie type of $\g$ and the divisibility of $\ell$
by $m$.   Fix a simple Lie algebra $\g$ of rank
$r$ and a positive integer $\ell$. Let $\la=\sum_i a_i\la_i$ be a
dominant weight of $\g$ written as an $\N$-linear combination of
fundamental weights $\la_i$.  To determine if $\la\in C_\ell(\g)$,
we compute:
$$\lan\la+\rho,\vartheta_j\ra=\lan\rho,\vartheta_j\ra+\sum_i^r a_i\lan\la_i,\vartheta_j\ra$$
where $j=0$ or $1$ depending on if $m\mid\ell$ or not.  Setting
$L_i^{(j)}=\lan\la_i,\vartheta_j\ra$ we see that the condition that
$\la\in C_\ell(\g)$ becomes:
$$\sum_i^k a_i L_i^{(j)}\leq\ell-\lan\rho,\vartheta_j\ra-1.$$  Since
$a_i,L_i^{(j)}\in\N$ we have:
\begin{lemma}
The cardinality of $C_\ell(\g)$ is the number of partitions of all
natural numbers $n$, $0\leq n\leq\ell-\lan\rho,\vartheta_j\ra-1$ into
parts from the size $r=\rm{rank}(\g)$ multiset
$\mathcal{S}(\g,\ell_m)=[L_i^{(j)}]_i^r$.
\end{lemma}

 So it remains only to compute the numbers $\lan\rho,\vartheta_j\ra$
and $L_i^{(j)}$ (with $j=0,1$) for each Lie algebra $\g$ and integer
$\ell>\lan\rho,\vartheta_j\ra$ and to apply standard combinatorics to
count the number of partitions into parts in
$\mathcal{S}(\g,\ell_m)$. The first task is easily accomplished with
the help of the book \cite{Bou}. Table \ref{rank} lists the results
of these computations, where $\ell_0:=\min\{\ell:\ell\geq\lan\rho,\vartheta_j\ra+1\}$ is the minimal non-degenerate
value of $\ell$.

\begin{table}\caption{$\mathcal{C}(\g,q,\ell)$ Data}\label{rank}
\begin{tabular}{*{3}{|c}|}
\hline
$X_r$ & $\mathcal{S}(\g,\ell_m)$ & $\ell_0$\\
\hline\hline
 $A_r$ & $[1,\ldots,1]$  & $r+1$
\\ \hline  $B_r$, $\ell$ odd& $[1,2,\ldots,2]$ & $2r+1$
\\ \hline  $B_r$, $\ell$ even& $[2,2,4,\ldots,4]$  & $4r-2$
\\ \hline  $C_r$, $\ell$ odd& $[1,2,\ldots,2]$  & $2r+1$
\\ \hline  $C_r$, $\ell$ even& $[2,\ldots,2]$  & $2r+2$
\\ \hline  $D_r$& $[1,1,1,2,\ldots,2]$  & $2r-2$
\\ \hline  $E_6$& $[1,1,2,2,2,3]$ & $12$
\\ \hline  $E_7$& $[1,2,2,2,3,3,4]$ & $18$
\\ \hline  $E_8$& $[2,2,3,3,4,4,5,6]$ & $30$
\\ \hline  $F_4$, $\ell$ even& $[2,4,4,6]$ & $18$
\\ \hline  $F_4$, $\ell$ odd& $[2,2,3,4]$ & $13$
\\ \hline  $G_2$, $3 \mid \ell$& $[3,6]$ & $12$
\\ \hline  $G_2$, $3\nmid \ell$& $[2,3]$ & $7$\\
\hline
\end{tabular}
\end{table}

Let $P_{\mathcal{T}}(n)$ denote the number of partitions of $n$ into
parts in a multiset $\mathcal{T}$, and
$P_{\mathcal{T}}[s]=\sum_{n=0}^sP_{\mathcal{T}}(n)$ the number of
partitions of all integers $0\leq n\leq s$ into parts from the
multiset $\mathcal{T}$.  Any standard reference on generating
functions (see e.g. \cite{Stan}) will provide enough details about
generating functions to prove the following:
\begin{lemma}\label{gflemma}
The number $P_{\mathcal{T}}(n)$ of partitions of $n$ into parts from
the multiset $\mathcal{T}$ has generating function:
$$\prod_{t\in\mathcal{T}}\frac{1}{1-x^t}=\sum_{n=0}^\infty P_{\mathcal{T}}(n)x^n,$$
while the number $P_{\mathcal{T}}[s]$ of partitions of all $n\in\N$
with $0\leq n\leq s$ into parts from the multiset $\mathcal{T}$ has
generating function:
$$\frac{1}{1-x}\prod_{t\in\mathcal{T}}\frac{1}{1-x^t}=\sum_{s=0}^\infty P_{\mathcal{T}}[s]x^s.$$
\end{lemma}
 Applying this lemma to the sets $\mathcal{S}(\g,\ell_m)$ given in
Table \ref{rank} we obtain:
\begin{theorem}
Define
$$F_{\g,\ell_m}(x)=\frac{1}{1-x}\prod_{k\in\mathcal{S}(\g,\ell_m)}\frac{1}{1-x^k}.$$
 Then the rank $|C_\ell(\g)|$ of the pre-modular category
$\mathcal{C}(\g,q,\ell)$ is
 the coefficient of  $$x^{\ell-\ell_0+\ell_m}$$ in the Taylor series expansion
 of $F_{\g,\ell_m}(x)$.
\end{theorem}
\begin{proof}
It is clear from Lemma \ref{gflemma} that the coefficients of
generating function $F_{\g,\ell_m}(x)$ counts the appropriate
partitions.  The coefficient of $x$ that gives the rank for a
specific $\ell$ is shifted by the minimal non-degenerate $\ell_0$,
which corresponds to the $x^0=1$ term if $m\mid\ell$ and to the
$x^1=x$ term of $m\nmid\ell$, hence the correction by $x^{\ell_m}$.
With this normalization only the coefficients of those powers of $x$
divisible (resp. ~indivisible)
 by $m$ give ranks corresponding to $\ell$ divisible
 (resp. ~indivisible) by $m$.
\end{proof}
We illustrate the application of this theorem with some examples.
\begin{example}
Let $\g$ be of type $G_2$. \begin{enumerate}
\item[(a)] Let $\ell=27$.
Then $\ell_m=0$ and $\ell_0=12$.  So the rank of
$\mathcal{C}(\g(G_2),q,27)$ is given by the $(27-12+0)=15$th
coefficient of
$$\frac{1}{(1-x)(1-x^3)(1-x^6)}=(1+x+x^2)(1+2x^3+4x^6+6x^9+9x^{12}+12x^{15}+\cdots)$$
so $|C_{27}(\g(G_2))|=12$.
\item[(b)] Let $\ell=14$.  Then $\ell_m=1$ and $\ell_0=7$.  So
$|C_{14}(\g(G_2))|$ is the $(14-7+1)$th coefficient of
$$\frac{1}{(1-x)(1-x^2)(1-x^3)}=1+x+2x^2+3x^3+4x^4+5x^5+7x^6+8x^7+10x^8\cdots$$
so the rank of $\mathcal{C}(\g(G_2),q,14)$ is $10$.
\end{enumerate}
\end{example}

\section{Examples}
We provide examples of two pre-modular categories, one of which is
modular and unitary, while the other is not modular but has a
(non-unitary) modular subcategory.  We only give enough information
to discuss the modularity and unitarity of the category.
\subsection{Type $\Z(A_1)$ at $\ell=5$}\label{ex1}
The following MTC is obtained from
$\mathcal{C}(\mathfrak{sl}_2,5,e^{\pi i/5})$ by taking the
subcategory of modules with integer highest weights. There are two
simple objects $\1$,and $X_1$ satisfying fusion rules: $X_1\otimes
X_1=\1\oplus X_1$ and $\1\otimes X_i=X_i$.
The $S$-matrix is $S=\begin{pmatrix} 1 & \frac{1+\sqrt{5}}{2}\\
\frac{1+\sqrt{5}}{2} & -1 \\
\end{pmatrix}$
and the twists: $\theta_0=1$, $\theta_1=e^{4\pi i/5}$.  It is clear
that $\det(S)\not=0$, and it follows from \cite{WenzlCstar} that the
category is unitary (notice that the categorical dimensions are both
positive).
\subsection{Type $B_2$ at $9$th Roots of Unity}
Consider the pre-modular categories
$\mathcal{C}(\mathfrak{so}_5,9,e^{j\pi i/9})$ with  $\gcd(18,j)=1$.
There are 12 inequivalent isomorphism classes of simple objects. The
simple iso-classes of objects are labelled by
$(\lambda_1,\lambda_2)\in\frac{1}{2}(\N^2)$ with
$\lambda_1\geq\lambda_2$.  The twist coefficients for $X_\la$ is
$q^{\lan \lambda+2\rho,\lambda\ra}$ where the form is twice the
usual Euclidean form.  The obstruction to modularity mentioned in
the proof of Theorem \ref{bmod} is labelled by
$\gamma:=\frac{1}{2}(5,5)$  The categorical dimension function is:
$$d_\la:=\frac{[2(\la_1+\la_2+2)][2(\la_1-\la_2+1)][2\la_1+3][2\la_2+1]}{[4][3][2][1]}.$$
One checks that the simple object $X_\gamma$ is indeed the cause of
the singularity of the $S$-matrix, that is,
$S_{\gamma,\lambda}=d_\gamma d_\la$ for all $\lambda$.  Thus this
category is not modular by Brugui\`eres' criterion, Theorem \ref{brug}.

Now let us consider the subcategory of
$\mathcal{C}(\mathfrak{so}_5,9,e^{j\pi i/9})$ with  $\gcd(18,j)=1$
generated by the simple objects labelled by integer weights:
$$\{(0,0),(1,0),(2,0),(1,1),(2,1),(2,2)\}.$$ The braiding and twists
from the full category restrict, so the entries of the $S$-matrix
are computed from Formula (\ref{sformula}). Taking the ordering of
simple objects above, we denote the categorical dimensions by $d_i$
$0\leq i\leq 5$. The fusion matrix corresponding to $(1,0)$ is:
$$N_1:=\begin{pmatrix}
0&1&0&0&0&0\\1&0&1&1&0&0\\
0&1&0&0&1&0\\0&1&0&1&1&0
\\0&0&1&1&1&1\\0&0&0&0&1&1
\end{pmatrix}.$$
It is not hard to show that $N_1$ determines the other five fusion
matrices by observing that $N_1$ has six distinct eigenvalues and
the fusion matrices commute. There are a total of six categories
corresponding to the six possible values of $q$. To describe the
$S$-matrices we let $\alpha$ be a primitive $18$th root of unity,
and set $r_1=-\alpha-\alpha^2+\alpha^5$,
$r_2=\alpha+\alpha^2-\alpha^4$ and $r_3=\alpha^4-\alpha^5$.  Then we
get the following $S$-matrices (for the 6 choices of $\alpha$):
$$\begin{pmatrix}
1&r_2& r_3&1&-1&r_1
\\r_2&1&1&r_1&-r_3&1
\\r_3&1&1&r_2&-r_1&1
\\1&r_1&r_2&1&-1
&r_3
\\-1&-r_3&-r_1&-1&1&-r_2\\r_1&1&1&r_3&-r_2&1
\end{pmatrix}.$$
One checks that $\det(S)\not=0$ for any $\alpha$, so these
categories \emph{are} modular. A bit of Galois theory shows that
there are only three distinct $S$ for the six choices of $\alpha$.
Notice that it is already clear that the first column of $S$ is
never positive, since both $1$ and $-1$ appear regardless of the
choice of $\alpha$.  So none of these categories is unitary.

\bibliographystyle{amsalpha}

\end{document}